\documentclass[12pt]{amsart}
\usepackage{amscd,amssymb}
\usepackage[arrow,matrix]{xy}

\theoremstyle{plain}
\newtheorem{thm}[subsection]{Theorem}
\newtheorem{lem}[subsection]{Lemma}
\newtheorem{prop}[subsection]{Proposition}
\newtheorem{cor}[subsection]{Corollary}

\theoremstyle{definition}
\newtheorem{rk}[subsection]{Remark}
\newtheorem{definition}[subsection]{Definition}
\newtheorem{ex}[subsection]{Example}

\numberwithin{equation}{section}
\newcommand{\F}{{\mathcal F}}
\newcommand{\A}{{\mathcal A}}

\newcommand{\B}{{\mathcal B}}

\newcommand{\LL}{{\mathcal L}}
\newcommand{\G}{{\mathcal G}}
\newcommand{\V}{{\mathcal V}}

\newcommand{\al}{{\alpha}}
\newcommand{\be}{{\beta}}

\newcommand{\Z}{\mathbb{Z}}
\newcommand{\Q}{\mathbb{Q}}
\newcommand{\R}{\mathcal{R}}
\newcommand{\C}{\mathbb{C}}

\newcommand{\PP}{\mathbb{P}}
\newcommand{\HH}{\mathbb{H}}

\newcommand{\T}{\mathbb{T}}

\DeclareMathOperator{\Hom}{Hom}

\begin{document}

\title [On admissible rank one local systems]
{ On admissible rank one local systems }

\author[Alexandru Dimca]{Alexandru Dimca }
\address{  Laboratoire J.A. Dieudonn\'e, UMR du CNRS 6621,
                 Universit\'e de Nice-Sophia-Antipolis,
                 Parc Valrose,
                 06108 Nice Cedex 02,
                 FRANCE.}
\email {dimca@math.unice.fr}

\subjclass[2000]{Primary 14C21, 14F99, 32S22 ; Secondary 14E05, 14H50.}

\keywords{local system, constructible sheaf, twisted cohomology, characteristic variety, resonance variety}

\begin{abstract} A rank one local system $\LL$ on a smooth complex algebraic variety $M$ is 1-admissible if the dimension of the first cohomology group $H^1(M,\LL)$
can be computed from the cohomology algebra $H^*(M,\C)$ in degrees $\leq 2$. Under the assumption that $M$ is 1-formal, we show that all local systems, except finitely many, on a non-translated irreducible component $W$ of the first characteristic variety $\V_1(M)$
are 1-admissible, see Proposition 3.1.
The same result holds for local systems on a translated component $W$, but now $H^*(M,\C)$ should be replaced by $H^*(M_0,\C)$, where $M_0$ is a Zariski open subset obtained from $M$ by deleting some hypersurfaces determined by the translated component $W$, see Theorem 4.3.
One consequence of this result is that the local systems $\LL$ where the dimension of  $H^1(M,\LL)$
jumps along a given positive dimensional component of the characteristic variety $\V_1(M)$
have finite order, see Theorem 4.7. Using this, we show in Corollary 4.9 that 
$\dim H^1(M,\LL)= \dim H^1(M,\LL^{-1})$ for any rank one local system $\LL$ on a smooth complex algebraic variety $M$.

\end{abstract}

\maketitle

\section{Introduction } \label{s0}

Let $M$ be a connected finite CW-complex. If $M$ is 1-formal, then the first twisted
Betti number of $M$ with coefficients in $\LL$ may be computed from the cohomology ring of $M$ in low
degrees, for rank one complex local systems $\LL$ near the trivial local system, see
\cite{DPS1}, Theorem A, the Tangent Cone Theorem.

In this paper, assuming moreover that $M$ is a connected smooth quasi-projective variety,
our aim  is to show that (a version of) the above statement is true
globally, with finitely many exceptions. In such a situation the exponential mapping \eqref{e1}
sends the irreducible components $E$ of the first resonance variety $\R_1(M)$ of $M$ onto
the non translated irreducible components $W$ of the first characteristic variety $\V_1(M)$ of $M$.

For $\al \in E$, $\al \ne 0$ (resp. $\LL \in W$, $\LL\ne \C_M$), the dimension of the cohomology group
$H^1(H^*(M,\C), \alpha \wedge)$ (resp.  $H^1(M,\LL)$) is constant (resp. constant with finitely many exceptions where this dimension may possibly increase). The first result is that the 1-formality assumption implies the inequality
\begin{equation} 
\label{e0}
\dim H^1(M,\LL) \ge \dim H^1(H^*(M,\C), \alpha \wedge)
\end{equation}
obtained by Libgober and Yuzvinsky when $M$ is a hyperplane arrangement complement, see \cite{LY}, Proposition 4.2.

An 1-admissible local system is a system for which the equality in the inequality \eqref{e0} holds. Various
characterization of 1-admissible local systems $\LL$ on a non translated component $W$ of the first characteristic variety $\V_1(M)$ are given in Proposition
\ref{p1}. In particular, we show that all local systems, except finitely many, on a non-translated irreducible component $W$ are 1-admissible.

 The main novelty is the analysis of
local systems belonging to a positive dimensional translated component $W'$ of the
first characteristic variety of $M$, see the last section. Such local systems (at least generically) are not 1-admissible. However, for a generic local system in $W'$, an equality similar to \eqref{e0} holds
but now $H^*(M,\C)$ should be replaced by $H^*(M_0,\C)$, where $M_0$ is a Zariski open subset obtained from $M$ by deleting some hypersurfaces determined by the translated component $W'$, see Theorem \ref{t3}.

One consequence of this result is the fact that the local systems $\LL$ where the dimension of  $H^1(M,\LL)$
jumps along a given positive dimensional irreducible component of the characteristic variety $\V_1(M)$
are local system of finite order, see Theorem \ref{t4}. Using this, we show in Corollary \ref{c21} that one has
$$\dim H^1(M,\LL)= \dim H^1(M,\LL^{-1})$$
 for any rank one local system $\LL$ on a smooth complex algebraic variety $M$.
In this section the role played by the constructible sheaf point of view introduced in \cite{D3} is
essential.

\section{Admissible and 1-admissible local systems} \label{s1}

Let $M$ be a smooth, irreducible, quasi-projective complex variety
and let $\T(M)=\Hom(\pi_1(M),\C^*)$ be the character variety of 
$M$. This is an algebraic group whose identity irreducible component is an algebraic torus
$\T(M)_1 \simeq (\C^*)^{b_1(M)}$. Consider the exponential mapping
\begin{equation} 
\label{e1}
\exp :H^1(M,\C) \to H^1(M,\C^*)=\T(M)
\end{equation}
induced by the usual exponential function $\exp: \C \to \C^*$. Clearly $\exp(H^1(M,\C))=\T(M)_1$.

\begin{definition} \label{d1}
A local system $\LL \in \T(M)_1$ is 1-admissible if there is a cohomology class $\alpha \in H^1(M,\C)$ such that $\exp(\alpha)=\LL$ and
$$ \dim H^1(M,\LL)=\dim H^1(H^*(M,\C), \alpha \wedge).$$
\end{definition}

If $\LL=\C_M$, then we can take $\alpha=0$ and the equality of dimension in Definition
\ref{d1} is obvious. So in the sequal we consider only the case $\LL \ne \C_M$.

\begin{rk} \label{r1}
When $M$ is a hyperplane arrangement complement or, more generally, a hypersurface 
arrangement complement in some projective space $\PP^n$, one usually defines the notion of
{\it admissible} local system $\LL$ on $M$ in terms of some conditions on the residues of an associated
logarithmic connection $\nabla(\al)$ on a good compactification of $M$, see for instance \cite{ESV},\cite{STV}, \cite{DM}.
For such an admissible local system $\LL$ on $M$ one has
$$ \dim H^i(M,\LL)=\dim H^i(H^*(M,\C), \alpha \wedge)$$
for all $i$ in the hyperplane arrangement case and for $i=1$ in the hypersurface arrangement case. For the case of  hyperplane arrangement complements, see also \cite{F} and \cite{LY}.
It is clear that "admissible" implies "1-admissible", which is a simpler, but still rather interesting property as we see below.
\end{rk}

One has the following easy result.

\begin{lem} \label{l1}
Any local system $\LL \in \T(M)$ is 1-admissible if $\dim M=1$.
\end{lem}

\proof Note that in this case the integral homology group $H_1(M)$ is torsion free and hence
$\T(M)=\T(M)_1$. 
Since $\LL$ is not the trivial local system, clearly one has $H^0(M,\LL)=0$. 

If $M$ is compact, then by duality, see \cite{D1}, we get
$$H^2(M,\LL)=H^0(M,\LL^{\vee})=0$$
and hence 
$$\dim H^1(M,\LL)=b_1(M)-2=-\chi(M).$$ 
If $M$ is not compact, then $M$ is homotopically equivalent to an 1-dimensional CW-complex, and hence
$H^2(M,\LL)=0$. In this case we get $\dim H^1(M,\LL)=b_1(M)-1=-\chi(M)$.
It follows that in both cases one has
\begin{equation} \label{e2}
\dim H^1(M,\LL)=-\chi(M)=\dim H^1(H^*(M,\C), \alpha \wedge).
\end{equation}
Note also that if $\LL=\C_M$ only the choice $\alpha=0$ is good, while in the case
$\LL \ne \C_M$ any choice for $\alpha$ satisfying $\exp(\alpha)=\LL$ is valid.

\endproof

To go further, we need the characteristic and resonance varieties, whose definition is recalled below.

The {\em characteristic varieties}\/ of $M$ 
are the jumping loci for the cohomology of $M$, with 
coefficients in rank~$1$ local systems:
\begin{equation} \label{e3}
\V^i_k(M)=\{\rho \in \T(M) \mid \dim H^i(M, \LL_{\rho})\ge k\}.
\end{equation}
When $i=1$, we use the simpler notation $\V_k(M)=\V^1_k(M)$.

The {\em resonance varieties}\/ of $M$ 
are the jumping loci for the cohomology of the complex  $H^*(H^*(M,\C), \alpha \wedge)$, namely:
\begin{equation} \label{e4}
\R^i_k(M)=\{\alpha \in H^1(M,\C) \mid \dim H^i(H^*(M,\C), \alpha \wedge)\ge k\}.
\end{equation}
When $i=1$, we use the simpler notation $\R_k(M)=\R^1_k(M)$.

\begin{ex} \label{ex0}
Assume that $\dim M=1$ and $\chi(M)<0$. Then it follows from the equality \eqref{e2}
that
$$\V_1(M)=\T(M) \text { and } \R_1(M)=H^1(M,\C).$$
\end{ex}

The more precise relation between the resonance and characteristic varieties can be summarized as follows, see \cite{DPS1}.

\begin{thm} \label{t2} Assume that $M$ is 1-formal. Then the irreducible components $E$ of the resonance variety
$\R_1(M)$ are linear subspaces in $H^1(M,\C)$ and the exponential mapping \eqref{e1} sends these irreducible components $E$ onto the irreducible components $W$ of $\V_1(M)$ with $1 \in W$.
\end{thm}

\begin{rk} \label{r11} 

The fact that $M$ is 1-formal depends only on the fundamental group $\pi_1(M)$, see for details \cite{DPS1}.
The class of 1-formal varieties is large enough, as it contains all the projective smooth varieties
and any hypersurface complement in $\PP^n$, see \cite{DPS1}. In fact, if the Deligne mixed Hodge structure
on $H^1(M,\Q)$ is pure of weight 2, then the smooth quasi-projective variety $M$ is 1-formal, see \cite{M}. The converse implication is not true, since any smooth quasi-projective curve obtained by deleting $k>1$ points from a projective curve of genus $g>0$ is 1-formal, but the corresponding
mixed Hodge structure
on $H^1(M,\Q)$ is not pure. Several examples of smooth quasi-projective varieties with a pure Deligne mixed Hodge structure on $H^1(M,\Q)$ are given in \cite{DL}.
\end{rk}

\medskip

In the sequel we concentrate ourselves on the strictly positive dimensional irreducible components of the first characteristic variety $\V_1(M)$. They have the following rather explicit description, given by
Arapura \cite{A}, see also Theorem 3.6 in \cite{D3}.

\begin{thm} \label{t1}
Let $W$ be a $d$-dimensional irreducible component of the first characteristic variety $\V_1(M)$, with $d>0$.
Then there is a regular morphism $f:M \to S$ onto a smooth curve $S=S_W$ with $b_1(S)=d$ such that the generic fiber $F$ of $f$ is connected, and a torsion character $\rho \in \T(M)$ such that the composition 
$$\pi_1(F) \stackrel{i_{\sharp}} \longrightarrow \pi_1(M) \stackrel{\rho} \longrightarrow \C^*,$$
where $i:F \to M$ is the inclusion, is trivial and
$$W=\rho \cdot f^*(\T(S)).$$
In addition, $\dim W=-\chi(S_W)+e$, with $e=1$ if $S_W$ is affine and $e=2$ if $S_W$ is proper.
If  $\LL \in W$, then $\dim H^1(M,\LL) \geq -\chi(S_W)$ and equality holds for all such $\LL$
with finitely many exceptions when $1 \in W$.

\end{thm}

If $1 \in W$, we say that $W$ is a {\it non-translated component} and then one can take $\rho =1$. If $1 \notin W$, we say that $W$ is a {\it translated component}.

\medskip

The following result  was obtained in \cite{LY}, Proposition 4.2 in the case of hyperplane arrangement complements.

\begin{prop} \label{p10}
Assume that $M$ is 1-formal and $\exp(\alpha)=\LL$. Then
$$\dim H^1(M,\LL) \geq \dim H^1(H^*(M,\C), \alpha \wedge).  $$ 

\end{prop}

\proof

If $\al \notin \R_1(M)$ or if $\LL$ is the trivial local system, then there is nothing to prove.
Otherwise the result follows directly from Proposition 6.6 in \cite{DPS1}.

\endproof

\begin{rk} \label{r2} 
(i) When $M$ is 1-formal, then  if $\LL$ is 1-admissible and $H^1(M,\LL)\ne 0$, then the cohomology class $\al$ realizing the conditions in Definition \ref{d1} is necessarily in an irreducible component $E=T_1W$ of $\R_1(M)$,
such that $\LL$ belongs to the non-translated irreducible component $W=\exp(E)$. For all $\LL \in \T(M)$, except finitely many, this component $W$ is uniquely determined by $\LL$, see \cite{N}.

(ii) Again when $M$ is 1-formal, this also shows that all the local systems on a translated component $W$ of $\V_1(M)$, possibly except finitely many located at the intersections of $W$ with non-translated components, are not 1-admissible. Indeed, all the examples in \cite{S1} suggest that the local systems
situated at the intersection of two (or several) irreducible components of $\V_1(M)$ are not 1-admissible.

\end{rk}

\begin{rk} \label{r10}
Note that Proposition \ref{p10} implies in particular
$$\exp(\R_k(M)) \subset \V_k(M)$$
for all $k$.
Since the differential of $\exp$ at the origin is the identity, this implies
$$\R_k(M) \subset TC_1\V_k(M).$$
Since the other inclusion always hold, see Libgober \cite{L10}, it follows that
the inequality in Proposition \ref{p10} implies the equality
$$\R_k(M)= TC_1\V_k(M).$$
If $M$ is not 1-formal, then the tangent cone $TC_1\V_k(M)$ can be strictly contained in $\R_k(M)$, see for instance Examples 5.11 and  9.1 in \cite{DPS1}.
It follows that the assumption 1-formal is needed to infer the inequality in Proposition \ref{p10}.
In other words, one may have
$$\dim H^1(M,\LL) < \dim H^1(H^*(M,\C), \alpha \wedge) $$
for some varieties $M$, which shows that the last claim in Proposition 4.5 in \cite{DM} fails for $k=0$ and
a general quasi-projective variety $M=Z \setminus D$. For instance, if $M=M_g$ is the surface constructed in
Example 5.11 in \cite{DPS1}, one has
$$2g-2=\dim  H^1(M,\LL) < \dim H^1(H^*(M,\C), \alpha \wedge)=2g-1$$
for all $\LL \ne \C_M$ and $\al \ne 0$. So in this case the only 1-admissible loacal system is the trivial local system $\C_M$.

\end{rk}

\begin{cor} \label{c1}
If $M$ is 1-formal, then any $\LL \in \T(M)_1$ with $H^1(M,\LL)=0$ is 1-admissible. More precisely, if $\LL=\exp(\alpha)$, then
$H^1(H^*(M,\C), \alpha \wedge)= 0$. 
\end{cor}

\proof 
Assume that $\LL=\exp(\alpha)$ and $H^1(H^*(M,\C), \alpha \wedge)\ne 0$. Then Proposition \ref{p10} gives a contradiction.
\endproof

The following result says that $\al \in \exp^{-1}(\LL)$ which occurs in Definition \ref{d1}
cannot be arbitrary in general.
\begin{prop} \label{p05}
Assume that $\R_1(M)\ne H^1(M,\C)$. Then for any local system $\LL \in \T(M)_1$ there are infinitely many
$\al \in \exp^{-1}(\LL)$ such that $H^1(H^*(M,\C), \alpha \wedge)=0$.
\end{prop}

\proof

Since $\LL \in \T(M)_1$, there is a cohomology class $\al_0 \in H^1(M,\C)$ such that $\exp(\al_0)=\LL$.
Then $\exp^{-1}(\LL)=\al_0 + \ker \exp$. We have to show that the set $(\al_0+  \ker \exp) \setminus \R_1(M)$
is infinite. The result follows from the following.

\begin{lem} \label{l05}
Consider the lattice $L_n=(2 \pi i) \cdot \Z^n \subset \C^n$ for $n \ge 1$. Then, for any point $\al \in \C^n$ and any subset  $A \subset \al+L_n$ such that $(\al +L_n) \setminus A$ is finite,
the Zariski closure
of $A$ is $\C^n$.
\end{lem}

\proof It is enough to show that any polynomial $g \in \C[x_1,...,x_n]$ such that $(\al +L_n) \setminus Z(g)$ is finite,
where $Z(g)$ is the zero-set of $g$, satisfies $g=0$.

The case $n=1$ is obvious. Assume the property is established for $n-1 \ge 1$ and consider the projection
$p:\C^n \to \C^{n-1}$, $(x_1,...,x_n)\mapsto (x_1,...,x_{n-1})$. Let $q=p|Z(g):Z(g) \to \C^{n-1}$.
It follows that $q(Z(g))$ contains a subset of $p(\al)+L_{n-1}$ with a finite complement, and the induction hypothesis implies that $q$ is a dominant mapping, i.e. the Zariski closure of $q(Z(g))$ is $\C^{n-1}$.
If $g\ne 0$, then $Z(g)$ is purely $(n-1)$-dimensional, and hence the generic fibers of $q$ are $0$-dimensional.
In other words, it exists a non-zero polynomial $h \in \C[x_1,...,x_{n-1}]$ such that $\dim q^{-1}(y)>0$ implies $h(y)=0$.
On the other hand, for any $y_0=p(\al)+v$ where $v \in L_{n-1}$, the fiber $q^{-1}(y_0)$ contains infinitely many points of the form
$$\al+v+2\pi i s e_n$$
with $s \in \Z$ and $e_n=(0,...,0,1)$. It follows that $\dim q^{-1}(y_0)>0$ and hence $h(y_0)=0$. The induction hypothesis implies that $h=0$,
a contradiction. This ends the proof of this Lemma and hence the proof of Proposition \ref{p05}.

\endproof

\medskip

In view of Corollary \ref{c1}, we consider in the sequal only local systems $\LL \ne \C_M$ such that
$\LL \in \V_1(M) \cap \T(M)_1$.

\section{Non-translated components and 1-admissible local systems} \label{s2}

Let $M$ be a smooth, quasi-projective complex variety. Let $W$ be an irreducible component
of $\V_1(M)$ such that $1 \in W$ and $\dim W>0$. Let $f:M \to S$ be the morphism onto a curve described in Theorem \ref{t1}, such that $W=f^*(\T(S))$. Note that  $\F:=R^0f_*(\C_M)=\C_S$ (since the generic fiber of $f$ is connected) and set $\G:=R^1f_*(\C_M)$.

\begin{prop} \label{p1} 

If $M$ is 1-formal, then the following three conditions on a local system $\LL=f^{-1}\LL' \in W$,are equivalent.

\noindent (i) $\LL$ is 1-admissible;

\noindent (ii) $\dim H^1(M,\LL)=\min _{\LL_1 \in W}\dim  H^1(M,\LL_1)$. (This minimum is called the {\rm generic dimension} of $H^1(M,\LL)$ along $W$.)

\noindent (iii) the natural morphism $f^*:H^1(S,\LL') \to H^1(M,\LL)$ is an isomorphism.

The condition

\noindent (iv) $H^0(S,\G \otimes \LL')=0$

implies the  condition $(iii)$ and they are equivalent when $S$ is affine.
Moreover, all these conditions are fulfilled by all $\LL \in W$ except finitely many.

\end{prop}

\proof

The equivalence of $(i)$ and $(ii)$ follows from Proposition \ref{p10} combined with the fact that
$(ii)$ holds for all local systems $\LL \in W$ except finitely many.

\medskip

For the definition of the morphism $f^*:H^1(S,\LL') \to H^1(M,\LL)$, see \cite{D1}, p. 54 and the references given there. The implication $(ii) \Rightarrow (iii)$ follows from the exact sequence
\begin{equation} \label{e5}
0 \to H^1(S,\LL') \to H^1(M, \LL ) \to H^0(S,\G\otimes \LL') 
\end{equation}
where the first morphism is precisely $f^*$ and the last morphism is surjective when
$S$ is affine or $\LL' \in \T(S)$ is generic, see Prop.4.3 in \cite{D3}.
It follows from Lemma \ref{l1} that the dimension of $H^1(S,\LL')$ is constant for $\LL'$ non-trivial.
Proposition 4.5 in \cite{D3}
gives the generic vanishing of the group $H^0(S,\G\otimes \LL') $. It follows that the minimal value for
$\dim H^1(M,\LL)$ is precisely $\dim H^1(S,\LL')$, and in such a case the monomorphism $f^*$ becomes an isomorphism.

\medskip

Conversely, assume that $(iii)$ holds.
Let $d=\dim W=b_1(S)$. Since the generic fiber of $f$ is connected, it follows that
\begin{equation} \label{e6}
f^*: H^1(S,\C) \to H^1(M, \C) 
\end{equation}
is injective. Let $\LL'=\exp(\omega)$ and note that then $\LL=\exp(\al)$, where $\al=f^*(\omega)$.
Using now the injectivity \eqref{e6} and Lemma \ref{l1}, it follows that
\begin{equation} \label{e7}
\dim H^1(S,\LL')=\dim \frac{\{\be \in H^1(S,\C)~~|~~\omega \wedge \be=0\}}{\C\cdot \omega}=\dim \frac{\{\gamma \in E~~|~~\al \wedge \gamma=0\}}{\C\cdot \al}
\end{equation}
where $E=f^*(H^1(S,\C))$ is a $d$-dimensional vector subspace in $H^1(M,\C)$.
In fact, it follows from Theorem \ref{t2} that $E$ is the irreducible component of $\R_1(M)$ corresponding to 
the irreducible component $W$ of $\V_1(M)$.

To show that $(i)$ holds, it is enough to show that
\begin{equation} \label{e8}
\{\gamma \in E~~|~~\al \wedge \gamma=0\}=\{\delta \in H^1(M,\C)~~|~~\al \wedge \delta=0\}.
\end{equation} 
Note that $\al \in \R_s(M)$ exactly when 
$$\dim \{\delta \in H^1(M,\C)~~|~~\al \wedge \delta=0\} \ge s+1.$$
Using \cite{DPS1}, it follows that $R_s(M)=\cup_iR^i$, where the union is over all the irreducible components $R^i$
of $\R_1(M)$ such that $\dim R^i>s+p(i)$, with $p(i)=0$ if the corresponding curve $S_i$ is not compact and $p(i)=1$
when the corresponding curve $S_i$ is compact.

Case 1. If $S$ is not compact, then clearly $\al \in (E \setminus 0) \subset (\R_{d-1}(M) \setminus R_d(M))$.
It follows that
$$\dim \{\delta \in H^1(M,\C)~~|~~\al \wedge \delta=0\} =d=\dim E$$
hence we get the equality \eqref{e8} in this case.

Case 2. If $S$ is compact, then clearly $\al \in (E \setminus 0) \subset (\R_{d-2}(M) \setminus R_{d-1}(M))$.
 It follows that
$$\dim \{\delta \in H^1(M,\C)~~|~~\al \wedge \delta=0\} =d-1=\dim \{\gamma \in E~~|~~\al \wedge \gamma=0\}    $$
hence we get the equality \eqref{e8} in this case as well.

The last claim follows directly from Proposition 6.6 in \cite{DPS1}.

\endproof

\section{Translated components and 1-admissible local systems} \label{s3}

Consider now the case of a translated component $W=\rho \cdot f^*(\T(S))$ and recall the notation from Theorem \ref{t1}. Let $\LL_0$ be the local system corresponding to $\rho$. We assume that $1 \notin W$ and this implies that the singular support $\Sigma (\F)$ of the constructible sheaf $\F=R^0f_*(\LL_0)$ is non-empty see Corollary 5.9 in \cite{D3}
(and coincides with the set of points $s \in S$ such that the stalk $\F_s$ is trivial, see Lemma 4.2 in  \cite{D3}).
We set as above $\G=R^1f_*(\LL_0)$ and recall the exact sequence
\begin{equation} \label{e9}
0 \to H^1(S,\F \otimes \LL') \to H^1(M, \LL_0 \otimes \LL ) \to H^0(S,\G\otimes \LL') 
\end{equation}
where $\LL=f^{-1}\LL'$ and the last morphism is surjective when
$S$ is affine or $\LL' \in \T(S)$ is generic, see Proposition 4.3 in \cite{D3}.
Moreover, one has

\medskip

(A) $ H^0(S,\G\otimes \LL') =0$ except for finitely many $\LL' \in \T(S)$, see Proposition 4.5 in \cite{D3}, and

\medskip

(B) $\F=Rj_*j^{-1}\F$, where $S_0=S \setminus \Sigma (\F)$ and $j:S_0 \to S$ is the inclusion, see \cite{D4}.

The proof of this last claim goes like that. It is known that a point $c \in S$ is in $ \Sigma(\F)$ if and only if for a small disc $D_c$ centered
at $c$, the restriction of the local system $\LL_{\rho}$ to the associated tube $T(F_c)=f^{-1}(D_c)$ about the fiber $F_c$
is non-trivial. Let $T(F_c)'=T(F_c) \setminus F_c$ and note that the inclusion $i: T(F_c)' \to T(F_c)$
induces an epimorphism at the level of fundamental groups.

Hence, if $c \in \Sigma(\F)$, then $\LL_{\rho}|T(F_c)' $ is a non-trivial rank one local system. In particular
$$H^0(T(F_c)',\LL_{\rho})=0.$$
If we apply the Leray spectral sequence to the locally trivial fibration
$$F \to T(F_c)' \to D_c'$$
where $D_c'=D_c \setminus \{c\}$, we get
$$H^0(D_c',\F)=H^0(T(F_c)',\LL_{\rho})=0.$$
It follows that $\F|D_c'$ is a non-trivial rank one local system. Hence $H^0(D_c',\F)=H^1(D_c',\F)=0$,
which proves the isomorphism $\F=Rj_*j^{-1}\F$.

\medskip

We deduce from (B) that the following more general isomorphism

\medskip

(C) $\F\otimes \LL'=Rj_*j^{-1}(\F \otimes \LL')$ for any $\LL' \in \T(S)$. In particular
\begin{equation} \label{e10}
H^1(S,\F \otimes \LL') =\HH^1(S, Rj_*j^{-1}(\F \otimes \LL'))=H^1(S_0,j^{-1}(\F \otimes \LL')),
\end{equation}
where the last isomorphism comes from Leray Theorem, see \cite{D1}, p.33.

Let $M_0=M \setminus f^{-1}(\Sigma (\F))$ and denote by $f_0:M_0 \to S_0$ the surjective morphism induced by $f$.

\begin{lem} \label{l2}
$$\LL_0|M_0 \simeq f_0^{-1}(\F|S_0).$$

\end{lem}

\proof For any local system $\LL_1 \in \T(M)$, there is a canonical adjunction morphism
$$a:f^{-1}f_*\LL_1 \to \LL_1$$
see \cite{KS}, (2.3.4), p. 91. In fact, in our situation, $f$ and $f_0$ are open mappings, so for any point $x \in M_0$ and $B_x$ a small open neighbourhood of $x$, one sees that the restriction morphism
$$a(B_x):\LL_0(f^{-1}(f(B_x)) \to \LL_0(B_x)=\C$$
is an isomorphism. Indeed, by Lemma 4.2 in  \cite{D3}, note that $s \in S_0$ if and only if the restriction
$\LL_0|T(F_s)$ is trivial, where $T(F_s)$ is a small open tube $f^{-1}(D_s)$ about the fiber $F_s=f^{-1}(s)$,
with $D_s$ a small disc centered at $s \in S$.

\endproof

The same proof yields the following more general result.

\begin{cor} \label{c2}
For any $\LL' \in \T(S)$ one has
$$(\LL_0 \otimes \LL)|M_0 \simeq f_0^{-1}((\F \otimes \LL')|S_0)$$
with $\LL=f^{-1}(\LL')$.

\end{cor}

For all local systems $\LL' \in \T(S)$ except finitely many, the exact sequence \eqref{e9} and the equality \eqref{e10} yield
\begin{equation} \label{e11}
H^1(M, \LL_0 \otimes \LL ) \simeq H^1(S,\F \otimes \LL')\simeq H^1(S_0, \LL'')
\end{equation}
where $\LL''=j^{-1}(\F \otimes \LL')=(\F \otimes \LL')|S_0$ is a rank one local system on $S_0$.

Note that the curve $S$ in Theorem \ref{t1} satisfies $\chi(S)\leq 0$ and hence $\chi(S_0)=\chi(S)- |\Sigma(\F)|<0$. It follows by Prop.1.7, Section V in Arapura \cite{A} that $W_0=f_0^*(\T(S_0)$ is an irreducible component in the characteristic variety
$\V_1(M_0)$ such that $1 \in W_0$ and $\dim W_0=b_1(S_0) \ge 2$.

With this notation, our main result is the following.

\begin{thm} \label{t3} Assume that $M$ is a smooth quasi-projective irreducible complex variety.
Let $W=\rho \cdot f^*(\T(S))$ be a translated $d$-dimensional irreducible component of the first characteristic variety $\V_1(M)$, with $d>0$.
Let $\LL_0$ be the rank one local system on $M$ corresponding to $\rho$, $\F=R^0f_*\LL_0$ and $\Sigma (\F)$ the singular support of $\F$. Set $S_0=S \setminus \Sigma (\F)$ and $M_0=f^{-1}(S_0)$. Assume moreover that
$M$ and $M_0$ are 1-formal.

Then there is a non-translated irreducible component $W_0$ of $\V_1(M_0)$, such that
$W \subset W_0$ under the obvious inclusion $\T(M) \to \T(M_0)$.
In particular, for any local system  $\LL_1 \in W$, except finitely many, there is a 1-form $\al(\LL_1)  \in H^1(M,\C)$ such that $\exp(\al(\LL_1))=\LL_1$ and $\dim H^1(H^*(M_0,\C), \al_0 (\LL_1) \wedge)=\dim H^1(M, \LL_1)$, where $\al_0(\LL_1)=\iota^*(\al(\LL_1))$, $\iota:M_0 \to M$ being the inclusion.

\end{thm}

\proof
With the above notation, apply Proposition \ref{p1} to the restriction $f_0:M_0 \to S_0$ and to the associated component $W_0$. We set $\LL_1=\LL_0 \otimes \LL$ and use \eqref{e11} and Corollary \ref{c2} to get
$$
\dim H^1(M, \LL_1 )= \dim H^1(S,\F \otimes \LL')= \dim H^1(S_0, \LL'')=$$
$$= \dim H^1(M_0, \LL_1|M_0)= \dim
H^1(H^*(M_0,\C), \al_0 (\LL_1) \wedge).$$
The key point here is that Proposition \ref{p1} holds for all local systems $\LL \in W$ except finitely many,
and not just for a generic local system in the sense of Zariski topology on $\T(M)$.

\endproof
\begin{rk} \label{r3}
One situation when clearly $M$ and $M_0$ are 1-formal is the following.
When $M$ is a hypersurface arrangement complement $M(\A)$ in some $\PP^n$, one can view $M_0$ as a new hypersurface arrangement complement $M(\B)$, where $\B$
is obtained from $\A$ by adding some additional components $H_W$ corresponding to the fibers in $f^{-1}(\Sigma(\F))$. Conversely, $\A$ is obtained from $\B$ by deleting the hypersurfaces in $H_W$.
So any translated component $W$ in $\V_1(M)$ corresponds to a non-translated component $W_0$ in a richer arrangement $\B=\A \cup H_W$. Even if $\A$ is a hyperplane arrangement, we see no reason why the richer arrangement $\B$ should contain only hyperplanes.
\end{rk} 

\begin{rk} \label{r4}
The dimension of $W_0$ is exactly the generic dimension of $H^1(M, \LL)$ for $\LL \in W$ plus one.
Indeed, one has $\dim W_0=-\chi(S_0)+1$, since $S_0$ is clearly non-compact, see \cite{D3}, Thm. 3.6.(i).
Moreover, the generic dimension of $H^1(M_0,\LL_0 \otimes \LL)$ is $-\chi(S_0)$, see \cite{D3}, Thm. 3.6.(iv).
On the other hand, since the generic dimension is realized outside a finite number of points on $W_0$, it follows that the generic dimension of $H^1(M_0,\LL_0 \otimes \LL)$ coincides to the generic dimension of  $H^1(M,\LL_0 \otimes \LL)$.

\end{rk} 

\begin{ex} \label{ex1} This is a basic example discovered by A. Suciu, see Example 4.1 in \cite{S1}, the so called deleted $B_3$-arrangement. Consider the line arrangement in $\PP^2$ given by the equation
$$xyz(x-y)(x-z)(y-z)(x-y-z)(x-y+z)=0.$$
Then there is a 1-dimensional translated component $W$. In this case the new hypersurface $H_W$ is the line
$x+y-z=0$, and $M_0$ is exactly the complement of the $B_3$-arrangement.

The characteristic variety $\V_1(M_0)$ has a 2-dimensional component $W_0$ denoted by $\Gamma$ in
Example 3.3 in \cite{S1}. In the notation of loc. cit. one has
$$\Gamma=\{(t,s,(st)^{-2},s,t,(st)^{-1},s^2,(st)^{-1})~~|~~(s,t)\in (\C^*)^2\}.$$
An easy computation shows that $W$ corresponds to the translated 1-dimensional torus
inside $W_0=\Gamma$ given by $st=-1$.

\end{ex}

The proof of Theorem \ref{t3} also gives the following result, which follows in the compact case from Simpson's work \cite{Si} and in the non-proper case from Budur's recent paper \cite{Bu}. 
In both cases one should also use in addition a result in \cite{DPS2}, saying that two irreducible components
of $\V_1(M)$ intersect at most at finitely many points, all of them torsion points in $\T(M)$.
Our proof below is
much simpler and purely topological.

\begin{thm} \label{t4} Assume that $M$ is a smooth quasi-projective irreducible complex variety.
Let $W$ be a  $d$-dimensional irreducible component of the first characteristic variety $\V_1(M)$, with $d>0$.
Let $\LL \in W$ be a rank one local system on $M$ such that 
$$\dim H^1(M,\LL) >\min _{\LL_1 \in W}\dim  H^1(M,\LL_1).$$
Then $\LL$ is a torsion point of the algebraic group $\T(M)$.
\end{thm}

\proof

Consider first the case when $W$ is a non-translated component. Then there is a morphism $f:M \to S$ and a local system $\LL' \in \T(S)$ such that $\LL=f^{-1}(\LL')$. The exact sequence \eqref{e5} implies that 
$$H^0(S,\G \otimes \LL')\ne 0$$
where we use the same notation as in \eqref{e5}.
On the other hand, Proposition 4.5 in \cite{D3} implies that
\begin{equation} \label{e20}
H^0(S_0,\G \otimes \LL')\ne 0
\end{equation}
where $S_0=S \setminus \Sigma(\G)$. Choose a finite set of generators $\gamma_1,...,\gamma_m$ for the group
$\pi_1(S_0)$. Then the condition \eqref{e20} implies that for each $j=1,...,m$, the monodromy $\lambda_j$ of
the local system $\LL'$ along the path $\gamma_j$ is the inverse of one of the eigenvalues of the monodromy operator $T_j$ of the geometric local system $\G|S_0$ along the path $\gamma_j$. Since the geometric local system $\G|S_0$ comes from an algebraic morphism, the Monodromy Theorem, see for instance \cite{De0},
implies that all the eigenvalues of any monodromy operator $T_j$ are roots of unity. Hence the same is true for all $\lambda_j$, which shows that $\LL'|S_0$ is a torsion point in $\T(S_0)$. The inclusion
$\T(S) \to \T(S_0)$ shows that $\LL'$ is a torsion point, and hence the same holds for $\LL=f^{-1}(\LL')$.

\medskip

Consider now the case when $W$ is a translated component. It follows from the proof of Theorem \ref{t3}
(for this part we do not need the 1-formality assumptions) that there is a Zariski open subset $M_0 \subset M$
and a non-translated component $W_0$ in $\V_1(M_0)$ such that $W \subset W_0$ under the natural inclusion
$\T(M) \subset \T(M_0)$ given by $\LL \mapsto \LL|M_0$.
Moreover, in this case $\LL=\LL_0 \otimes f^{-1}(\LL')$ for some morphism $f:M \to S$ and a local system $\LL' \in \T(S)$. The subset $M_0$ depends on the torsion local system $\LL_0$ and on $f$, but not on $\LL'$.
The equality
$$\dim H^1(M, \LL)= \dim H^1(M_0, \LL|M_0)$$
for all such local systems $\LL$ obtained by varying $\LL'$, and the fact that along any 
component the jumps in dimension occur only at finitely many points, recall Theorem \ref{t1} for non-translated components and use Corollary 5.9 in \cite{D3} in the translated case, implies that
$\LL$ is a jumping point for dimension along $W$ if and only if $\LL|M_0$ is a jumping point for dimension along $W_0$. The first part of this proof shows that in such a case $\LL|M_0$ is a torsion point, and  
the inclusion
$\T(M) \to \T(M_0)$ shows that the same holds for $\LL$.

\endproof

\begin{cor} \label{c3}
Let $M$ be a smooth, quasi-projective complex variety which is 1-formal. Let $W$ be an irreducible component
of $\V_1(M)$ such that $1 \in W$ and $\dim W>0$. If $\LL \in W$ is not 1-admissible, then $\LL$ is a torsion point in $\T(M)$.

\end{cor}

Note that the algebraic group of characters $\T(M)$ has a complex conjugation involution,
denoted by $\LL \mapsto \overline \LL$ and satisfying
$$\dim H^k(M,\LL)=\dim H^k(M,\overline \LL)$$
for all $k$. This follows simply by noting that the complex of finite dimensional $\C$-vector spaces used to compute the twisted cohomology $H^*(M,\LL)$ (resp. $H^*(M,\overline \LL)$) comes from a complex of real vector spaces
such that the corresponding differentials $d^k(\LL)$ and $d^k(\overline \LL)$ are complex conjugate
for all $k$.

This remark  and the above Theorem \ref{t4} have the following consequence.

\begin{cor} \label{c21}
Let $M$ be a quasi-projective manifold. Then, for any rank one local system $\LL \in \T(M)$
one has $\dim H^1(M,\LL)= \dim H^1(M,\LL^{-1})$.
\end{cor}

\proof We can assume that either $\LL$ or $\LL^{-1}$ is in $\V_1(M)$. Since the situation is symmetric,
assume that $\LL \in \V_1(M)$. If $\LL$ belongs to a strictly positive dimensional component $W$ of
$\V_1(M)$, it follows from the description of such components given in \cite{D3}, Corollary 5.8, that $W^{-1} \subset \V_1(M)$. Moreover, the generic dimension of $H^1(M,\LL)$ along $W$ and along $W^{-1}$
coincide by \cite{D3}, Corollary 5.9. If $\LL$ or $\LL^{-1}$ is a jumping point for this dimension, it follows from Theorem \ref{t4} that both $\LL$ and $\LL^{-1}$ are torsion points in the group $\T(M)$. Hence $\LL^{-1}=\overline \LL$ and in this latter case the claim follows from the above remark.

On the other hand, if $\LL$ is an isolated point of $\V_1(M)$, it follows from \cite{A}
that $\LL$ corresponds to a unitary character, and hence  again $\LL^{-1}=\overline \LL$
and we conclude as above.

\endproof

\end{document}